\documentclass[preprint,12pt,nopreprintline]{elsarticle}
\usepackage{amsmath,amssymb,amsthm}
\usepackage{mleftright}
\mleftright
\usepackage{graphicx}
\graphicspath{{./}{./figures/}}
\usepackage{booktabs}
\usepackage{multirow}
\usepackage{float}
\usepackage{caption}
\usepackage{subcaption}
\usepackage{algorithm}
\usepackage{algpseudocode}
\usepackage{tikz}
\usetikzlibrary{positioning,arrows.meta,fit}
\usepackage{hyperref}
\usepackage{geometry}
\theoremstyle{definition}
\newtheorem{lemma}{Lemma}
\newtheorem{proposition}{Proposition}

\allowdisplaybreaks

\begin{document}
\begin{frontmatter}

\title{A multi-objective mixed integer linear programming model for supply chain planning of 3D printing}

\author{Amirreza Talebi}

\affiliation{organization={Department of Integrated Systems Engineering},
            addressline={The Ohio State University},
            city={Columbus},
            postcode={43210},
            state={OH},
            country={USA}}

\begin{abstract}
Three-dimensional (3D) printing makes customized parts in small batches, but 3D printers are expensive and slow, so the grouping, orientation, and sequencing of parts decide both delivery times and printer use. We study the allocation of parts with due dates to jobs on identical powder bed 3D printers, where the orientation of a part sets its height and the plate area it occupies. We formulate a bi-objective mixed integer linear programming model that minimizes the weighted earliness and tardiness of the parts and the unused plate area of the jobs, and we show that the optimal timing of a fixed job sequence is an isotonic regression. The Pareto front is computed with the epsilon-constraint method, and a non-dominated sorting genetic algorithm built on this timing result handles larger instances. On random instances, choosing the orientation lowers the earliness and tardiness by 63\% on average compared with the orientation given in the data. For 15 to 40 parts, the genetic algorithm keeps 98\% to 99\% of the hypervolume of the best known fronts in less than a minute per run on average, whereas a time-limited epsilon-constraint method keeps about half of it with 30 and 40 parts.
\end{abstract}


\end{frontmatter}

\section{Introduction}\label{sec:intro}

Three-dimensional (3D) printing, also called additive manufacturing (AM), builds a part layer by layer from a digital model, without molds or tooling. It makes small batches of customized parts economical, and it has moved from the production of prototypes to the production of final parts in industrial, medical, and aerospace applications \cite{berman2012,niaki2018management,gebler2014global}. 3D printers are still expensive and slow compared with mass production lines, so the planning of their production decides how fast orders reach customers and how well the costly machines are used.

In powder bed printers, several parts are printed together in one job. The time of a job is a setup time plus a term that grows with the material volume of the parts and a term that grows with the job height, which is the height of the tallest part \cite{li2017production,kucukkoc2019milp}. A part can be placed in the build chamber in several orientations. An orientation that lowers the height of a part enlarges the area that the part occupies on the build plate \cite{oh2020impact,che2021machine}. Grouping parts into jobs, orienting them, and sequencing the jobs on the printers are therefore linked decisions.

We consider customers who want each part on its due date. A part that is finished early waits in storage, and a part that is finished late delays the customer, so the planner wants just-in-time delivery. The planner also wants each job to use as much of the build plate as possible, because every job costs printer time, material handling, and energy. These two goals conflict. More jobs let parts finish close to their own due dates, while fewer and fuller jobs use the plate better.

The contributions of this paper are as follows. First, we formulate a bi-objective mixed integer linear programming (MILP) model that combines job formation, part orientation, the sequencing of jobs on identical printers, and the earliness and tardiness of parts, with the unused plate area as the second objective. Second, we show that the unused area equals the plate area times the number of jobs minus the total base area of the parts, so the Pareto front is organized by the number of jobs, and we derive a bound on completion times that keeps an optimal solution and tightens the big-M constraints. Third, we show that, for a fixed sequence of jobs, the optimal completion times solve an isotonic regression problem that the pool adjacent violators algorithm (PAVA) solves exactly. We use this result in a non-dominated sorting genetic algorithm II (NSGA-II). Fourth, we compute exact Pareto fronts with the epsilon-constraint method, measure the value of orientation against three fixed orientation policies, study the sensitivity to the layer time, the plate size, and the number of printers, and compare NSGA-II with the epsilon-constraint method on instances with up to 40 parts. The code and data are provided with the paper.

Section~\ref{sec:lit-review} reviews the literature. Section~\ref{sec:problemdefinition} defines the problem, Section~\ref{sec:mathmodeling} presents the model, its properties, and the two solution methods, and Section~\ref{sec:numexp} reports the experiments. Section~\ref{sec:conclusion} concludes.

\section{Literature review}\label{sec:lit-review}

Planning for AM differs from classical batch scheduling because the parts of a job share the setup and the layer time, and the job height is set by its tallest part. \cite{pinto2024nesting} map the research on nesting and scheduling in AM with a bibliometric and systematic review. We review the studies closest to our problem.

\cite{li2017production} allocate parts to machines with different production times and costs per volume, times per unit height, setup times, and maximum supported areas and heights, so as to minimize the average production cost per volume. They solve a mathematical model with CPLEX and propose the best-fit and adapted best-fit heuristics. \cite{kucukkoc2019milp} develops MILP models that minimize the makespan for a single AM machine, identical parallel machines, and non-identical parallel machines. \cite{li2019dynamic} study dynamic order acceptance and scheduling for on-demand production in powder bed systems, with the goal of maximizing the average profit per unit time. \cite{kim2022part} group parts by material and schedule the builds on non-identical parallel machines with sequence-dependent setup times to minimize the makespan.

Due dates appear in several studies. \cite{chergui2018production} propose a heuristic based on the earliest due date (EDD) rule for nesting and scheduling parts on identical parallel AM machines when orders have due dates. \cite{fera2018modified} model the scheduling of a single selective laser melting machine to minimize time and cost and solve it with a genetic algorithm that works on a matrix encoding, tested on a production plan with 30 parts. \cite{10.1007/978-3-030-76307-7_8} minimize the total tardiness on non-identical selective laser melting machines that differ in dimensions, speed, and cost, with parts that differ in size, release date, and due date, using a genetic algorithm. \cite{dvorak2018planning} print parts with unique configurations and deadlines on several machines and combine bin packing, two-dimensional nesting, and job shop scheduling in a constraint model. \cite{luzon2019job} choose job sizes and sequences when the printing process can fail and a failed job restarts from the beginning, with an application to dental printing. \cite{ransikarbum2017multi} assign parts to fused deposition modeling printers with a multi-objective model that balances the operating cost, the load of the printers, the total tardiness, and the number of unprinted parts, assuming a fixed part orientation.

Part orientation changes the height and the base area of a part \cite{leutenecker2016considering}. \cite{zhang2016two} combine orientation optimization with a two-dimensional placement of several parts. \cite{oh2020impact} compare a laying policy, which minimizes part heights, with a standing policy, which minimizes the projected base areas so that fewer jobs are needed, and find that the standing policy can give shorter completion times when about 40 parts are printed. \cite{che2021machine} select orientations while scheduling unrelated parallel batch machines to minimize the makespan, with a MILP model and a simulated annealing algorithm, and \cite{hu2022scheduling} add two-dimensional packing and release times to this setting and use an adaptive large neighborhood search for large instances.

These studies minimize the makespan, the tardiness, or a cost, mostly with one objective, and those that choose orientations do not consider earliness. To the best of our knowledge, no study combines the earliness and tardiness of parts, the use of the build plate, and the choice of part orientation in one bi-objective model. This paper fills this gap.

\section{Problem definition}\label{sec:problemdefinition}

A set of parts with given due dates has to be printed on identical powder bed 3D printers. A job is a set of parts printed together on one printer in one run. Each part is printed in exactly one job, each printer processes its jobs one after another, and all parts of a job are finished when the job is finished. Figure~\ref{fig:1} shows an allocation of parts to jobs and printers.

\begin{figure}[H]
\centering
\begin{tikzpicture}[
  part/.style={circle,draw,fill=blue!12,minimum size=7mm,inner sep=0pt,font=\small},
  job/.style={rectangle,draw,rounded corners,fill=gray!15,minimum width=19mm,minimum height=10mm,align=center,font=\small},
  prn/.style={rectangle,draw,fill=green!20,minimum width=17mm,minimum height=10mm,font=\small},
  arr/.style={-{Latex[length=2mm]}}]
\foreach \k in {1,...,9} \node[part] (p\k) at (0,{-0.75*(\k-1)}) {$P_{\k}$};
\node[draw,dashed,rounded corners,fit=(p1)(p9),inner sep=4pt] (pool) {};
\node[prn] (m1) at (3.4,-0.6) {Printer 1};
\node[prn] (m2) at (3.4,-3.0) {Printer 2};
\node[prn] (m3) at (3.4,-5.4) {Printer 3};
\foreach \q in {1,2,3} \draw[arr] (pool.east) -- (m\q.west);
\node[job, right=7mm of m1] (j11) {Job 1\\$P_1, P_3$};
\node[job, right=5mm of j11] (j12) {Job 2\\$P_4, P_8$};
\node[job, right=7mm of m2] (j21) {Job 1\\$P_6, P_7$};
\node[job, right=5mm of j21] (j22) {Job 2\\$P_2$};
\node[job, right=7mm of m3] (j31) {Job 1\\$P_5, P_9$};
\draw[arr] (m1) -- (j11); \draw[arr] (j11) -- (j12);
\draw[arr] (m2) -- (j21); \draw[arr] (j21) -- (j22);
\draw[arr] (m3) -- (j31);
\draw[arr] ([yshift=-9mm]j31.south west) -- ++(5.4,0) node[right,font=\small] {time};
\end{tikzpicture}
\caption{Allocation of parts to jobs on identical 3D printers and sequencing of the jobs on each printer (adapted from \cite{kucukkoc2019milp}).}
\label{fig:1}
\end{figure}

\begin{figure}[H]
\centering
\begin{subfigure}{0.49\textwidth}
\centering
\includegraphics[width=\textwidth]{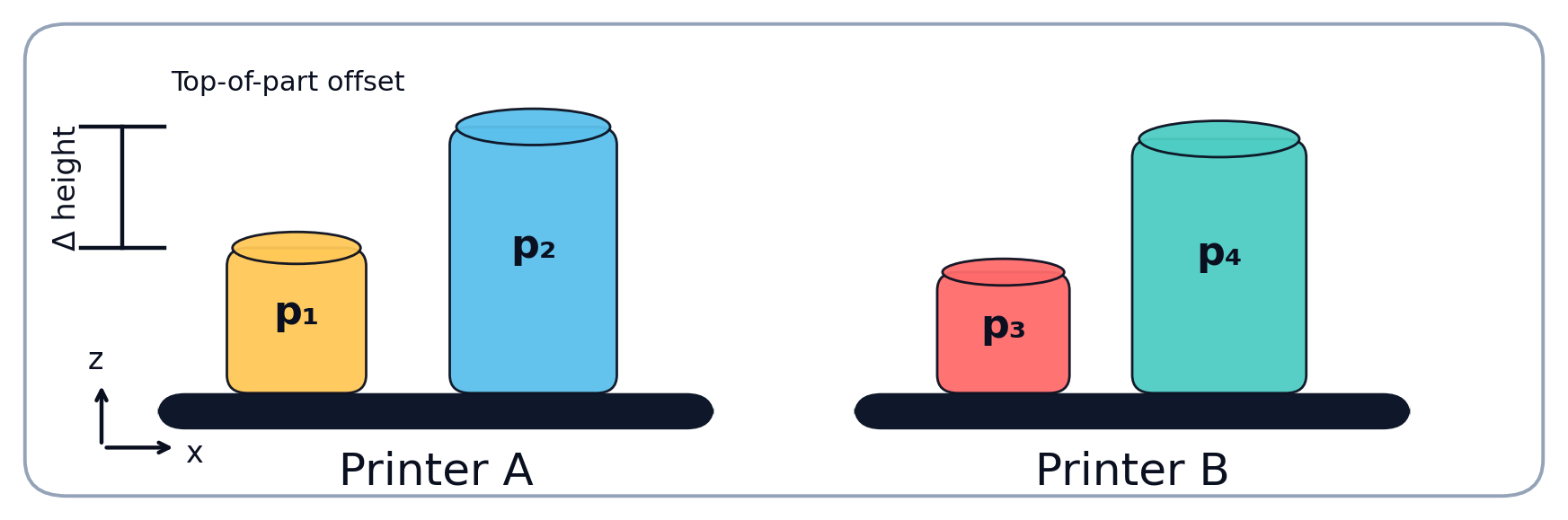}
\caption{Combination 1}
\end{subfigure}
\hfill
\begin{subfigure}{0.49\textwidth}
\centering
\includegraphics[width=\textwidth]{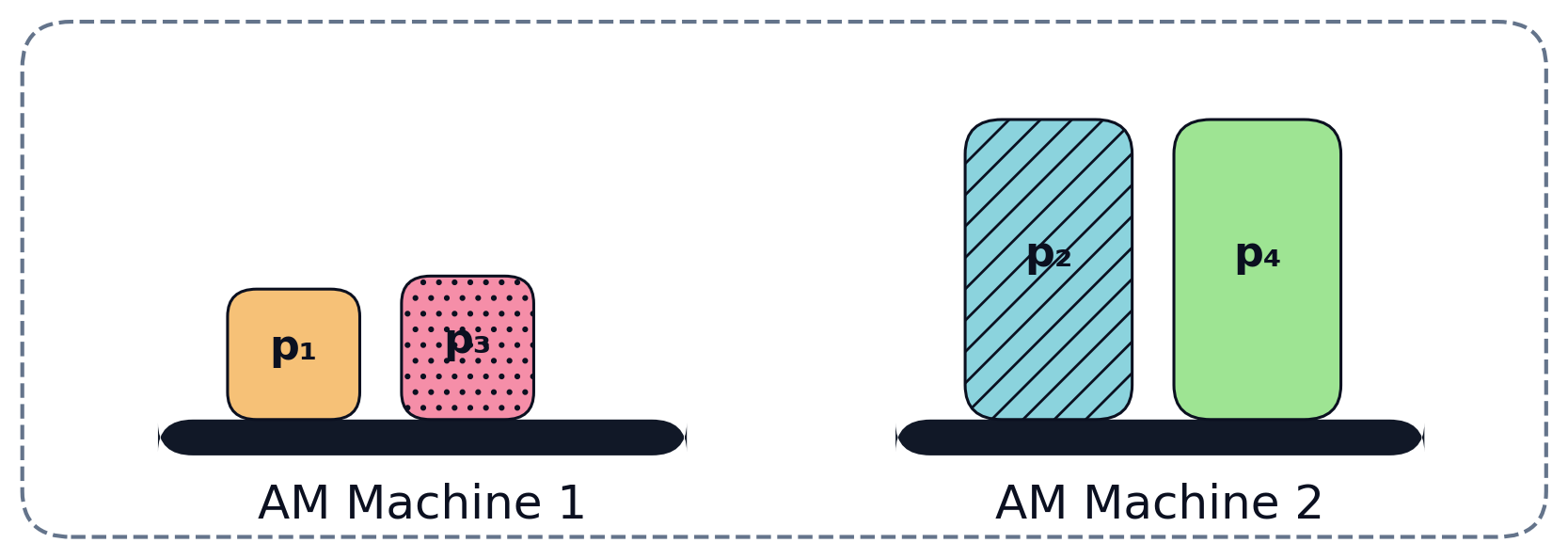}
\caption{Combination 2}
\end{subfigure}
\caption{Two ways to combine four parts into two jobs (adapted from \cite{oh2018production}). In combination 1, the tall parts $P_2$ and $P_4$ set the completion of both jobs. In combination 2, the short parts share a job and finish earlier.}
\label{fig:combinations}
\end{figure}

The processing time of a job is the sum of a setup time, a layer time proportional to the height of its tallest part, and a volume time proportional to the material volume of its parts \cite{li2017production,kucukkoc2019milp}. Because the tallest part sets the layer time of the whole job, parts of similar height and similar due date should be printed together. Figure~\ref{fig:combinations} illustrates this point. In the first combination, the tall parts $P_2$ and $P_4$ set the completion of both jobs, and the short parts $P_1$ and $P_3$ wait for them. In the second combination, the two short parts share one job and finish early, and the two tall parts share the other job.

\begin{figure}[H]
\centering
\includegraphics[width=13cm]{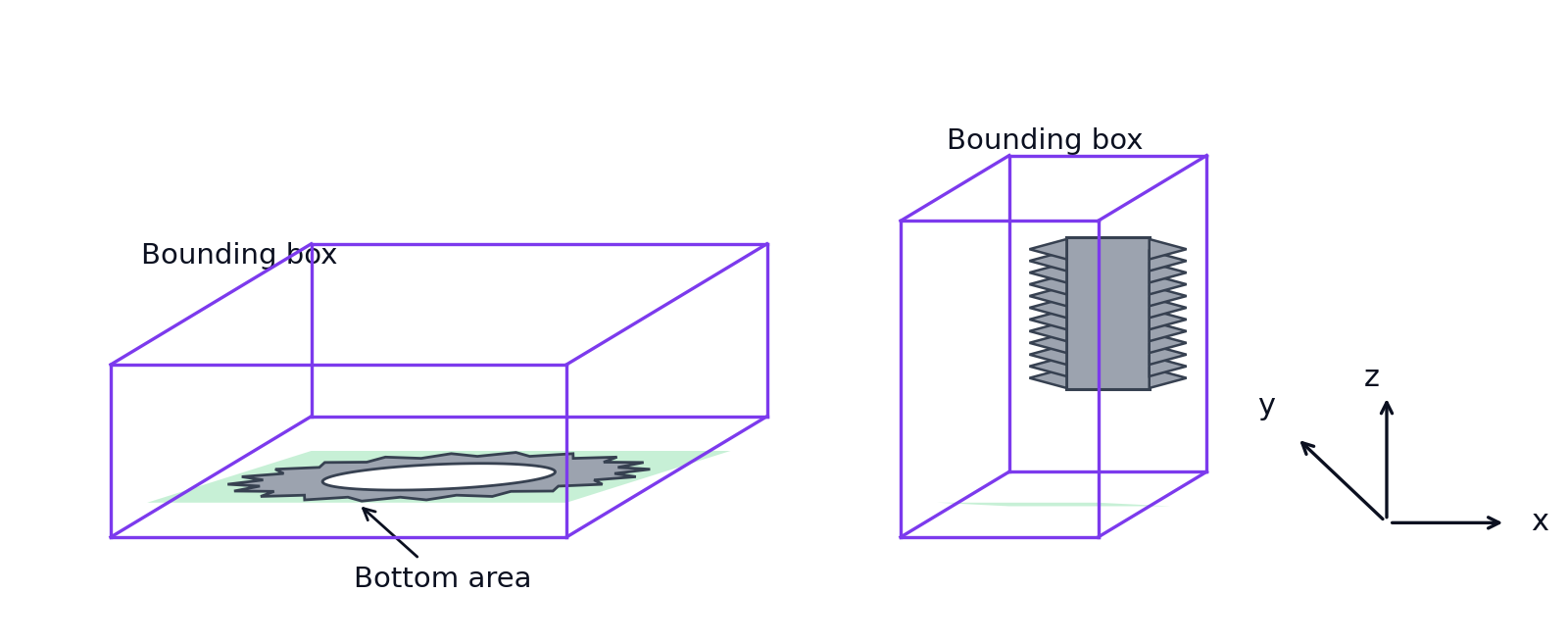}
\caption{Bounding box and base area of a part in two orientations (adapted from \cite{oh2018production}).}
\label{fig:4}
\end{figure}

We make the following assumptions.
\begin{enumerate}
\item Each part is represented by its bounding box with edges $l_i$, $w_i$, and $h_i$, and it is printed with one of these three edges vertical (Figure~\ref{fig:4}). A rotation about the vertical axis does not change the height or the base area of the part, so three orientations are enough.
\item The capacity of a job is the usable area of the build plate. The sum of the base areas of the parts in a job may not exceed this area, as in \cite{li2017production,kucukkoc2019milp}. The usable area is smaller than the plate area, to leave room for spacing between parts and for the losses of two-dimensional nesting.
\item An orientation is admissible if the part height does not exceed the chamber height and the base of the part fits on the plate.
\item All parts are made of the same material and are available at time zero, the printers are identical, and a job cannot be interrupted.
\item A printer may stay idle between two jobs, so that parts are not finished too early.
\end{enumerate}

The first objective is the weighted sum of the earliness and tardiness of the parts. The second objective is the unused usable area of the build plate, summed over the jobs, which measures how well the printers are used.

\section{Mathematical model}\label{sec:mathmodeling}

\subsection{Notation}

Table~\ref{tab:notation} lists the sets, parameters, and decision variables. Orientation $o=1$ is the orientation given in the data, and orientations $o=2$ and $o=3$ put the edges $l_i$ and $w_i$ upright:
\begin{equation}\label{eq:orient}
(h_{i1},a_{i1})=(h_i,l_iw_i), (h_{i2},a_{i2})=(l_i,w_ih_i), (h_{i3},a_{i3})=(w_i,l_ih_i).
\end{equation}
For every orientation, $h_{io}a_{io}=l_iw_ih_i$, so a lower part has a larger base.

\begin{table}[H]
\centering
\caption{Notation.}
\label{tab:notation}
\small
\begin{tabular}{lp{13.2cm}}
\toprule
Symbol & Description \\
\midrule
\multicolumn{2}{l}{Sets and indices}\\
$\mathcal{I}$, $i$ & Parts, $n=|\mathcal{I}|$ \\
$\mathcal{M}$, $m$ & Identical 3D printers \\
$\mathcal{J}$, $j$ & Job positions on each printer, $\mathcal{J}=\{1,\dots,J\}$; job $(j,m)$ is the $j$th job of printer $m$ \\
$\mathcal{O}_i$, $o$ & Admissible orientations of part $i$, $\mathcal{O}_i\subseteq\{1,2,3\}$ \\
\midrule
\multicolumn{2}{l}{Parameters}\\
$l_i,w_i,h_i$ & Edges of the bounding box of part $i$ (mm) \\
$h_{io},a_{io}$ & Height (mm) and base area (mm$^2$) of part $i$ in orientation $o$, see \eqref{eq:orient} \\
$v_i$ & Material volume of part $i$ (mm$^3$) \\
$d_i$ & Due date of part $i$ (h) \\
$\bar{A},\bar{H}$ & Usable plate area (mm$^2$) and chamber height (mm) \\
$\tau^S,\tau^H,\tau^V$ & Setup time per job (h), time per mm of job height (h/mm), time per mm$^3$ of material (h/mm$^3$) \\
$\alpha,\beta$ & Weights of earliness and tardiness (per hour) \\
$\bar{C}$ & Upper bound on the completion times of jobs (h), see Lemma~\ref{lem:horizon} \\
\midrule
\multicolumn{2}{l}{Decision variables}\\
$x_{ijmo}\in\{0,1\}$ & 1 if part $i$ is printed in job $(j,m)$ in orientation $o$ \\
$y_{jm}\in\{0,1\}$ & 1 if job $(j,m)$ is used \\
$H_{jm}\geq 0$ & Height of job $(j,m)$ (mm) \\
$P_{jm}\geq 0$ & Processing time of job $(j,m)$ (h) \\
$C_{jm}\geq 0$ & Completion time of job $(j,m)$ (h) \\
$E_i,T_i\geq 0$ & Earliness and tardiness of part $i$ (h) \\
\bottomrule
\end{tabular}
\end{table}

\subsection{Bi-objective MILP model}

The model reads
\begin{align}
\min\quad &f_1=\sum_{i\in\mathcal{I}}\left(\alpha E_i+\beta T_i\right), \label{eqn:f1}\\
\min\quad &f_2=\sum_{j\in\mathcal{J}}\sum_{m\in\mathcal{M}}\Bigl(\bar{A}y_{jm}-\sum_{i\in\mathcal{I}}\sum_{o\in\mathcal{O}_i}a_{io}x_{ijmo}\Bigr), \label{eqn:f2}\\
\text{s.t.}\quad &\sum_{j\in\mathcal{J}}\sum_{m\in\mathcal{M}}\sum_{o\in\mathcal{O}_i}x_{ijmo}=1, \forall i\in\mathcal{I}, \label{eqn:assign}\\
&\sum_{i\in\mathcal{I}}\sum_{o\in\mathcal{O}_i}a_{io}x_{ijmo}\leq\bar{A}y_{jm}, \forall j\in\mathcal{J},m\in\mathcal{M}, \label{eqn:area}\\
&y_{jm}\leq\sum_{i\in\mathcal{I}}\sum_{o\in\mathcal{O}_i}x_{ijmo}, \forall j\in\mathcal{J},m\in\mathcal{M}, \label{eqn:open}\\
&\sum_{o\in\mathcal{O}_i}h_{io}x_{ijmo}\leq H_{jm}\leq\bar{H}y_{jm}, \forall i\in\mathcal{I},j\in\mathcal{J},m\in\mathcal{M}, \label{eqn:height}\\
&P_{jm}=\tau^Sy_{jm}+\tau^HH_{jm}+\tau^V\sum_{i\in\mathcal{I}}\sum_{o\in\mathcal{O}_i}v_ix_{ijmo}, \forall j\in\mathcal{J},m\in\mathcal{M}, \label{eqn:proc}\\
&C_{1m}\geq P_{1m}, \forall m\in\mathcal{M}, \label{eqn:first}\\
&C_{jm}\geq C_{j-1,m}+P_{jm}, \forall j\in\mathcal{J}\setminus\{1\},m\in\mathcal{M}, \label{eqn:seq}\\
&y_{jm}\leq y_{j-1,m}, \forall j\in\mathcal{J}\setminus\{1\},m\in\mathcal{M}, \label{eqn:order}\\
&E_i\geq d_i\sum_{o\in\mathcal{O}_i}x_{ijmo}-C_{jm}, \forall i\in\mathcal{I},j\in\mathcal{J},m\in\mathcal{M}, \label{eqn:early}\\
&T_i\geq C_{jm}-d_i-(\bar{C}-d_i)\Bigl(1-\sum_{o\in\mathcal{O}_i}x_{ijmo}\Bigr), \forall i\in\mathcal{I},j\in\mathcal{J},m\in\mathcal{M}, \label{eqn:tardy}\\
&x_{ijmo},y_{jm}\in\{0,1\}, 0\leq C_{jm}\leq\bar{C}, H_{jm},P_{jm},E_i,T_i\geq 0, \forall i\in\mathcal{I},j\in\mathcal{J},m\in\mathcal{M},o\in\mathcal{O}_i. \label{eqn:domain}
\end{align}

Objective \eqref{eqn:f1} is the weighted earliness and tardiness of the parts. Early parts cause storage costs and late parts cause customer dissatisfaction or penalties, and the weights $\alpha$ and $\beta$ set their relative importance. Objective \eqref{eqn:f2} is the unused usable plate area of the used jobs. Constraints \eqref{eqn:assign} print each part once, in one job and one admissible orientation. Constraints \eqref{eqn:area} limit the base area of a job to the usable plate area and open the job when it holds a part, and constraints \eqref{eqn:open} close the jobs that hold no part, so that $\sum_{j,m}y_{jm}$ is the number of jobs. Constraints \eqref{eqn:height} set the job height to at least the height of each of its parts, and to zero for a closed job. Since $\sum_{o}h_{io}x_{ijmo}$ is the height of part $i$ when it is in job $(j,m)$ and zero otherwise, no big-M term is needed. Constraints \eqref{eqn:proc} give the processing time of a job. Constraints \eqref{eqn:first} and \eqref{eqn:seq} sequence the jobs of each printer, and they allow idle time between jobs. Constraints \eqref{eqn:order} use the job positions of a printer in order, which removes equivalent solutions with empty positions between used jobs. For a part $i$ in job $(j,m)$, constraints \eqref{eqn:early} and \eqref{eqn:tardy} give $E_i\geq d_i-C_{jm}$ and $T_i\geq C_{jm}-d_i$, and they are inactive otherwise, because $C_{jm}\geq 0$ and $C_{jm}\leq\bar{C}$. Since both objectives push $E_i$ and $T_i$ down, $E_i=(d_i-C_{jm})^+$ and $T_i=(C_{jm}-d_i)^+$ at an optimum.

The model is linear. The product of the orientation of a part and its assignment, which makes the area and height constraints nonlinear when the two decisions have separate variables, is avoided by indexing the assignment variable by the orientation. The completion time of a part, which is the product $C_{jm}\sum_o x_{ijmo}$, enters only through the big-M constraints \eqref{eqn:early} and \eqref{eqn:tardy}, whose constants are tightened by Lemma~\ref{lem:horizon} below.

In the implementation, we add three groups of valid inequalities. They are $x_{ijmo}\leq y_{jm}$, the bound $T_i\geq(p^{\min}_i-d_i)^+$, where $p^{\min}_i=\tau^S+\tau^H\min_{o\in\mathcal{O}_i}h_{io}+\tau^Vv_i$ is the shortest job that contains part $i$, and $x_{ijmo}=0$ for $m>i$, which labels the identical printers by the smallest index of the parts they print.

\subsection{Properties of the model}

\begin{proposition}\label{prop:f2}
For every feasible solution, let $N=\sum_{j,m}y_{jm}$ be the number of jobs and $o(i)$ the orientation of part $i$. Then
\begin{equation}\label{eq:f2form}
f_2=\bar{A}N-\sum_{i\in\mathcal{I}}a_{io(i)}=\bar{A}N-\sum_{i\in\mathcal{I}}\frac{l_iw_ih_i}{h_{io(i)}}\geq 0 .
\end{equation}
Hence $N\geq\bigl\lceil\sum_i\min_{o\in\mathcal{O}_i}a_{io}/\bar{A}\bigr\rceil$, and for a fixed number of jobs, $f_2$ decreases when parts are laid flatter.
\end{proposition}

\begin{proof}
By \eqref{eqn:assign}, $\sum_{j,m}\sum_{i,o}a_{io}x_{ijmo}=\sum_i\sum_o a_{io}\sum_{j,m}x_{ijmo}=\sum_i a_{io(i)}$, which gives the first equality, and $a_{io}=l_iw_ih_i/h_{io}$ gives the second. Summing \eqref{eqn:area} over the jobs gives $f_2\geq 0$ and $\bar{A}N\geq\sum_i a_{io(i)}\geq\sum_i\min_o a_{io}$, and $N$ is an integer.
\end{proof}

Proposition~\ref{prop:f2} explains the shape of the Pareto front. A new job adds $\bar{A}$ to $f_2$, while a change of orientation changes $f_2$ only by the difference of two base areas. The nondominated points therefore form groups with the same number of jobs. Laying parts flat lowers both the job heights, which helps $f_1$, and $f_2$. The conflict comes from the number of jobs, since more jobs let parts finish closer to their due dates, and from the plate capacity, since flat parts occupy more area.

The problem is also hard to solve. If all parts are cubes, every orientation has the same base area, and minimizing $f_2$ is the bin packing problem of packing the base areas into the fewest jobs of capacity $\bar{A}$, which is NP-hard \cite{bernhard2006bin}. This motivates the heuristic of Section~\ref{sec:nsga}.

For $\varepsilon\in\mathbb{R}\cup\{\infty\}$, let $P(\varepsilon)$ be the problem of minimizing $f_1$ subject to \eqref{eqn:assign} to \eqref{eqn:domain} and $f_2\leq\varepsilon$. The next lemma gives the constant $\bar{C}$.

\begin{lemma}\label{lem:horizon}
Let $d^{\max}=\max_i d_i$ and $\bar{p}_i=\tau^S+\tau^H\max_{o\in\mathcal{O}_i}h_{io}+\tau^Vv_i$. Suppose that $P(\varepsilon)$ without the bound $C_{jm}\leq\bar{C}$ has an optimal solution.
\begin{enumerate}
\item[(a)] It has an optimal solution with $C_{jm}\leq d^{\max}+\sum_{i\in\mathcal{I}}\bar{p}_i$ for all $j$ and $m$.
\item[(b)] If $\mathit{UB}$ is the value of $f_1$ at a feasible solution of $P(\varepsilon)$, it has an optimal solution with $C_{jm}\leq d^{\max}+\mathit{UB}/\beta$, $E_i\leq\mathit{UB}/\alpha$, and $T_i\leq\mathit{UB}/\beta$ for all $i$, $j$, and $m$.
\end{enumerate}
So $\bar{C}$ can be set to either bound without losing optimality.
\end{lemma}

\begin{proof}
Take an optimal solution. By \eqref{eqn:height} and \eqref{eqn:proc}, $H_{jm}=P_{jm}=0$ for every closed job, and by \eqref{eqn:order} the closed jobs of a printer come after its used jobs. (a) Let $D_{jm}=\sum_{j'\leq j}P_{j'm}$ and $C'_{jm}=\min\{C_{jm},d^{\max}+D_{jm}\}$. Then $C'_{1m}\geq P_{1m}$, and since $C_{jm}\geq C_{j-1,m}+P_{jm}$ and $d^{\max}+D_{jm}=d^{\max}+D_{j-1,m}+P_{jm}$, we get $C'_{jm}\geq\min\{C_{j-1,m},d^{\max}+D_{j-1,m}\}+P_{jm}=C'_{j-1,m}+P_{jm}$, so $C'$ is feasible. If $C'_{jm}<C_{jm}$, then $C'_{jm}=d^{\max}+D_{jm}\geq d_i$ for every part $i$ of the job, so these parts stay on time or tardy and their tardiness decreases. Hence $f_1$ does not increase, and $f_2$ does not change. Finally, each used job satisfies $P_{jm}\leq\sum_{i\text{ in job }(j,m)}\bar{p}_i$, so $C'_{jm}\leq d^{\max}+D_{jm}\leq d^{\max}+\sum_i\bar{p}_i$. (b) The optimal value is at most $\mathit{UB}$, so $\alpha E_i\leq\mathit{UB}$ and $\beta T_i\leq\mathit{UB}$ for every part. A used job $(j,m)$ holds some part $i$, so $C_{jm}\leq d_i+T_i\leq d^{\max}+\mathit{UB}/\beta$. The completion time of a closed job can be set to that of the job before it.
\end{proof}

The same bounds hold for the second stage of the epsilon-constraint method in Section~\ref{sec:eps}, whose optimal solutions have $f_1=f_1^*(\varepsilon)\leq\mathit{UB}$. Part (b) gives a much smaller constant than part (a) when a good solution is known. In the epsilon-constraint method of Section~\ref{sec:eps}, we take $\mathit{UB}$ from a solution with the smallest $f_2$, which is feasible for every $\varepsilon$ that is used.

The last property concerns the completion times. Fix the jobs, the orientations, and the sequence of the jobs on a printer, and let the jobs of the printer be $k=1,\dots,K$, with processing times $P_k$ and part sets $B_k$. Write $D_k=\sum_{k'\leq k}P_{k'}$, which is the completion time of job $k$ when the printer never waits, and $t_i=d_i-D_k$ for $i\in B_k$.

\begin{proposition}\label{prop:timing}
The completion times that minimize $\sum_k\sum_{i\in B_k}\bigl(\alpha(d_i-C_k)^++\beta(C_k-d_i)^+\bigr)$ subject to \eqref{eqn:first} and \eqref{eqn:seq} are $C_k=D_k+s_k$, where $s$ solves
\begin{equation}\label{eq:iso}
\min_{s}\sum_{k=1}^{K}g_k(s_k), \text{ s.t. } 0\leq s_1\leq s_2\leq\dots\leq s_K, \text{ with } g_k(s)=\sum_{i\in B_k}\bigl(\alpha(t_i-s)^++\beta(s-t_i)^+\bigr).
\end{equation}
PAVA solves \eqref{eq:iso} exactly. It keeps blocks of consecutive jobs. The value of a block $B$ is the $r$th smallest target $t_i$ of its parts, with $r=\lceil\alpha|B|/(\alpha+\beta)\rceil$. Jobs are added one at a time as new blocks, a block is merged with the block before it while the value of the earlier block is larger, and finally $s_k$ is the larger of zero and the value of the block that contains job $k$.
\end{proposition}

\begin{proof}
With $C_k=D_k+s_k$, constraint \eqref{eqn:first} becomes $s_1\geq 0$, and \eqref{eqn:seq} becomes $s_k\geq s_{k-1}$ because $D_k-D_{k-1}=P_k$. The cost of a part in job $k$ is $\alpha(t_i-s_k)^++\beta(s_k-t_i)^+$, which gives \eqref{eq:iso}. Each $g_k$ is convex and piecewise linear, so \eqref{eq:iso} minimizes a separable convex function over a chain, and PAVA is exact for this class \cite{best1990active,best2000minimizing}. For the parts of a block, the right derivative of their total cost at $s$ is $-\alpha|\{i:t_i>s\}|+\beta|\{i:t_i\leq s\}|$, which becomes nonnegative first at the $r$th smallest target, so this target minimizes the block cost. The lower bound $s_1\geq 0$ can be written as a job $0$ with $s_0\leq s_1$ and cost $W|s_0|$, where $W>(\alpha+\beta)n$. Then $s_0=0$ at the optimum, PAVA merges into the block of job 0 exactly the blocks with a negative value, and the result is the clipped solution.
\end{proof}

Proposition~\ref{prop:timing} extends the timing algorithm of \cite{garey1988one} for a fixed sequence of single jobs with symmetric weights to jobs that contain several parts and to unequal weights. We checked it against the linear program of the timing problem on 2000 random job sequences, and the largest difference in the optimal value was $2.3\times10^{-13}$.

\subsection{Epsilon-constraint method}\label{sec:eps}

A feasible solution $x$ dominates a feasible solution $x'$ if $f_k(x)\leq f_k(x')$ for $k=1,2$ with at least one strict inequality. A solution is Pareto optimal if no feasible solution dominates it, its objective vector is a nondominated point, and the set of nondominated points is the Pareto front \cite{hwang1980mathematical}. The epsilon-constraint method \cite{haimes1971bicriterion,mavrotas2009effective} finds nondominated points by solving $P(\varepsilon)$ for several values of $\varepsilon$. To avoid weakly Pareto optimal solutions, each solve has two stages. The first stage solves $P(\varepsilon)$ and gives $f_1^*(\varepsilon)$, and the second stage minimizes $f_2$ subject to the same constraints and $f_1\leq f_1^*(\varepsilon)$.

With integer edge lengths, $f_2$ takes integer values. The complete method starts with $\varepsilon=\infty$ and, after each nondominated point $(f_1^k,f_2^k)$, sets $\varepsilon=f_2^k-1$, until $P(\varepsilon)$ is infeasible. When every solve is optimal, the $k$th iteration returns the nondominated point with the $k$th largest $f_2$, so the method returns the whole Pareto front. For larger instances, we use the grid version, which solves the two lexicographic problems that give the extreme points and then $P(\varepsilon)$ for $\varepsilon$ on a grid of 11 equally spaced values between their $f_2$ values. In both versions, we first minimize $f_2$, time the solution optimally with Proposition~\ref{prop:timing}, use its $f_1$ as $\mathit{UB}$ in Lemma~\ref{lem:horizon}(b), and give this solution to the solver as a starting point. The completion times of every solution returned by the solver are also reset by Proposition~\ref{prop:timing} for its job sequences. This keeps $f_2$ and the value of an optimal solution, and it can lower $f_1$ when a solve stops at its time limit with slack in $E_i$ and $T_i$.

\subsection{NSGA-II}\label{sec:nsga}

NSGA-II \cite{deb2002fast} keeps a population of solutions, ranks them by fast non-dominated sorting, and keeps diversity by the crowding distance. We use it with the following encoding and decoder.

A chromosome has three parts. The first is a permutation $\pi$ of the parts, which is the order in which parts are packed into jobs. The second is an admissible orientation $o_i$ for each part. The third is a cut bit $c_i\in\{0,1\}$ for each part, where $c_i=1$ closes the current job after part $i$. Algorithm~\ref{alg:decode} decodes a chromosome. Parts are packed in the order $\pi$, and a new job is opened when the next part does not fit on the plate or when the previous part has a cut bit of one, so few cut bits give few full jobs and many cut bits give many small jobs. The jobs are sent, in the order they are formed, to the printer with the least total processing time among the printers with fewer than $J$ jobs, and the completion times on each printer are set by PAVA (Proposition~\ref{prop:timing}). A chromosome that needs more than $|\mathcal{M}|J$ jobs violates the job limit by the number of extra jobs, and the constrained domination rule of \cite{deb2002fast} ranks it after all feasible solutions.

\begin{algorithm}[H]
\caption{Decoder of a chromosome $(\pi,o,c)$.}
\label{alg:decode}
\begin{algorithmic}[1]
\State $\mathcal{B}\gets$ empty list of jobs, $B\gets\emptyset$
\For{$i$ in the order $\pi$}
  \If{$B\neq\emptyset$ and $\sum_{k\in B}a_{ko_k}+a_{io_i}>\bar{A}$} append $B$ to $\mathcal{B}$, $B\gets\emptyset$ \EndIf
  \State $B\gets B\cup\{i\}$
  \If{$c_i=1$} append $B$ to $\mathcal{B}$, $B\gets\emptyset$ \EndIf
\EndFor
\State append $B$ to $\mathcal{B}$ if $B\neq\emptyset$
\For{each job $B\in\mathcal{B}$ in order}
  \State $P_B\gets\tau^S+\tau^H\max_{i\in B}h_{io_i}+\tau^V\sum_{i\in B}v_i$
  \State send $B$ to the printer with the least load among the printers with fewer than $J$ jobs
\EndFor
\State time the jobs of each printer with PAVA and compute $f_1$, and set $f_2=\bar{A}|\mathcal{B}|-\sum_i a_{io_i}$
\end{algorithmic}
\end{algorithm}

The initial population mixes seeded and random chromosomes. One quarter uses the EDD order and the flattest admissible orientations with cut bits drawn at a random density, and it includes the two chromosomes with no cuts and with a cut after every part. One quarter perturbs the EDD order by adjacent swaps and draws each orientation as the flattest one or a random one, and the rest is random. Parents are chosen by binary tournaments on rank and crowding distance. With probability 0.9, a pair of parents produces two children by order crossover on $\pi$ and uniform crossover on $o$ and $c$. Each child then receives an insertion move on $\pi$ with probability 0.5, each orientation is reset with probability $1/n$ (to the flattest orientation or a random admissible one, with equal probability), and each cut bit is flipped with probability $1/n$. The parents and children are merged and the best half by rank and crowding distance survives \cite{deb2002fast}. The algorithm returns the nondominated feasible solutions of the final population.

\section{Numerical experiments}\label{sec:numexp}

\subsection{Instances and settings}\label{sec:setup}

Table~\ref{tab:params} lists the parameters of the instances. The build plate is 250 mm by 250 mm, of which 80\% is usable, so $\bar{A}=50{,}000$ mm$^2$. For a job of four flat parts, the setup time, the layer time, and the volume time are of similar size, which is the setting in which the orientation matters. Each part has one long edge and two short edges, and the three edges are put in random order to form $(l_i,w_i,h_i)$, so the orientation given in the data is random. The due dates are drawn around $\Lambda=\frac{1}{|\mathcal{M}|}\sum_i\bigl(\tau^Vv_i+(\tau^S+\tau^H\min_{o}h_{io})/4\bigr)$, which estimates the makespan of a schedule with four flat parts per job, so that some parts can be finished on time and others cannot. With $\alpha=\beta=1$, $f_1$ is the total earliness and tardiness in hours. All orientations of the generated parts are admissible on the 250 mm plate. Each printer has $J=\lceil n/(2|\mathcal{M}|)\rceil+1$ job positions, which allows two parts per job on average and one spare position.

\begin{table}[H]
\centering
\caption{Parameters of the instances.}
\label{tab:params}
\small
\begin{tabular}{ll}
\toprule
Parameter & Value \\
\midrule
Build plate, usable share, $\bar{A}$ & 250 mm by 250 mm, 0.8, 50,000 mm$^2$ \\
Chamber height $\bar{H}$ & 250 mm \\
Setup time $\tau^S$ & 2 h \\
Layer time $\tau^H$ & 0.07 h per mm of job height \\
Volume time $\tau^V$ & $10^{-5}$ h per mm$^3$ \\
Weights $\alpha$ and $\beta$ & 1 and 1 \\
Long edge & integer, uniform on 60 to 200 mm \\
Short edges (two) & integer, uniform on 15 to 80 mm \\
Material volume $v_i$ & $\rho_il_iw_ih_i$, with $\rho_i$ uniform on $[0.2,0.6]$ \\
Due date $d_i$ & uniform on $[0.6\Lambda,1.6\Lambda]$, rounded to 0.1 h \\
Job positions $J$ & $\lceil n/(2|\mathcal{M}|)\rceil+1$ \\
\bottomrule
\end{tabular}
\end{table}

Instance $r$ with $n$ parts and $|\mathcal{M}|$ printers uses the random seed $20261022+1000n+100|\mathcal{M}|+r$. The MILP is solved by HiGHS 1.15 \cite{huangfu2018parallelizing} through its Python interface, with one thread per solve, a relative optimality gap of $10^{-6}$, and a time limit of 600 s per solve unless stated otherwise. NSGA-II is written in Python, with a population of 200 and 500 generations. All runs were made on a laptop with an AMD Ryzen 7 7730U processor with 8 cores, with 12 runs in parallel, so the reported times are those of a loaded machine. The code, the instances, and all results are provided as supplementary material.

\subsection{Illustrative instance}\label{sec:illustrative}

\begin{figure}[H]
\centering
\includegraphics[width=\textwidth]{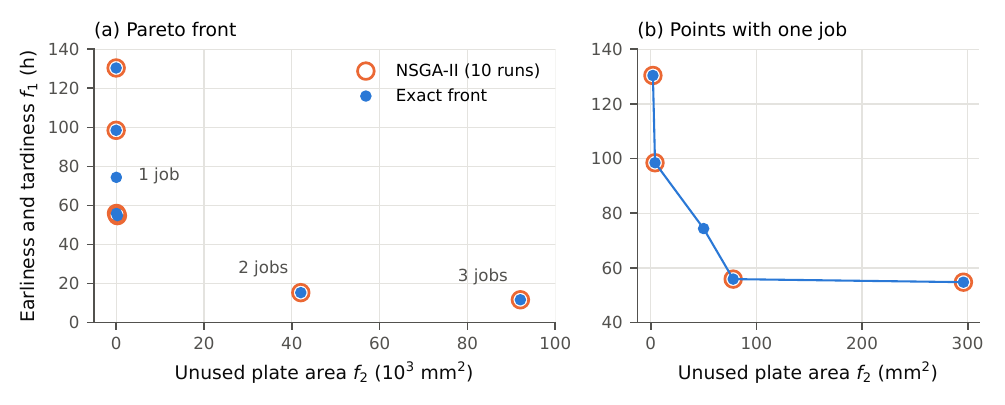}
\caption{Pareto front of the illustrative instance. Panel (b) enlarges the points with one job, which lie near $f_2=0$ in panel (a).}
\label{fig:front}
\end{figure}

We first solve instance 0 with 8 parts and 2 printers, whose data are given in Table~\ref{tab:example}. The complete epsilon-constraint method found the whole Pareto front, 7 nondominated points, with 8 two-stage solves in 71 s, and every solve was optimal. Table~\ref{tab:front} lists the points and Figure~\ref{fig:front} plots them.

\begin{table}[H]
\centering
\caption{Data of the illustrative instance (8 parts, 2 printers). The edges are in mm, the volume in cm$^3$, and the due date in h.}
\label{tab:example}
\small
\begin{tabular}{rrrrrr}
\toprule
Part & $l_i$ & $w_i$ & $h_i$ & $v_i$ & $d_i$ \\
\midrule
1 & 47 & 171 & 80 & 327.2 & 9.4 \\
2 & 150 & 38 & 41 & 129.1 & 6.2 \\
3 & 58 & 39 & 142 & 161.2 & 11.0 \\
4 & 18 & 22 & 199 & 45.5 & 11.8 \\
5 & 16 & 66 & 58 & 22.3 & 8.7 \\
6 & 22 & 78 & 109 & 51.9 & 6.3 \\
7 & 114 & 52 & 60 & 192.1 & 8.0 \\
8 & 64 & 44 & 99 & 60.3 & 14.8 \\
\bottomrule
\end{tabular}
\end{table}

\begin{table}[H]
\centering
\caption{Pareto front of the illustrative instance. The last column gives the number of the 10 NSGA-II runs that found the point.}
\label{tab:front}
\small
\begin{tabular}{rrrrr}
\toprule
Jobs & $f_1$ (h) & $f_2$ (mm$^2$) & Solve time (s) & NSGA-II runs \\
\midrule
3 & 11.72 & 92,050 & 65.5 & 10 \\
2 & 15.38 & 42,050 & 1.8 & 10 \\
1 & 54.81 & 296 & 0.5 & 4 \\
1 & 55.93 & 78 & 0.6 & 3 \\
1 & 74.41 & 50 & 0.8 & 0 \\
1 & 98.49 & 4 & 1.1 & 10 \\
1 & 130.41 & 2 & 0.9 & 3 \\
\bottomrule
\end{tabular}
\end{table}

The front has the structure that Proposition~\ref{prop:f2} predicts. The point with 3 jobs has the smallest earliness and tardiness, 11.72 h, and leaves 92,050 mm$^2$ of usable plate area unused. Using 2 jobs saves one plate, 50,000 mm$^2$, at the cost of 3.66 h. With a single job, all parts finish at the same time, and $f_1$ rises to at least 54.81 h. The five points with one job differ only in the orientations. Filling the plate up to 2 mm$^2$ of its usable area raises $f_1$ from 54.81 h to 130.41 h, because the plate can only be filled this tightly by standing tall parts upright, which raises the job height and the layer time. For a planner, the useful choices are the lowest point of each group, while the other points trade a negligible area for long delays.

\begin{figure}[H]
\centering
\includegraphics[width=\textwidth]{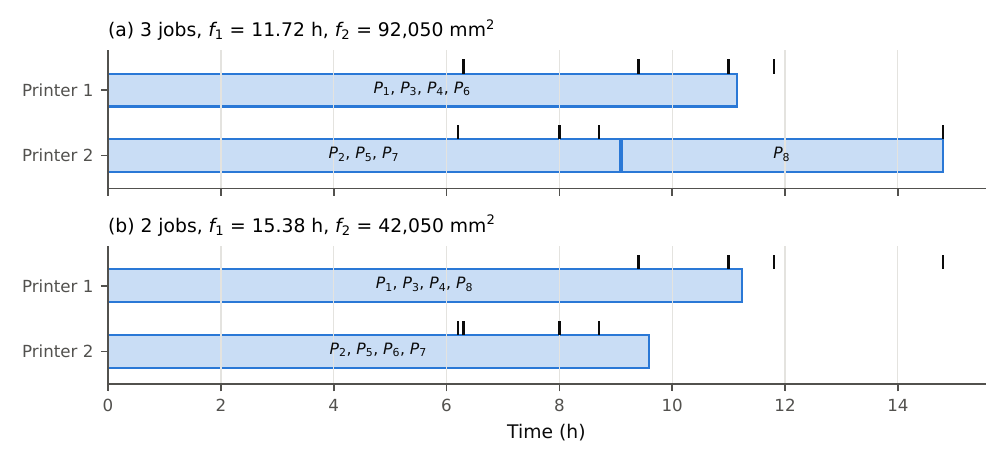}
\caption{Schedules of the Pareto points with 3 jobs (a) and 2 jobs (b). Each bar is a job, labeled with its parts, and the ticks above the bars mark the due dates of the parts of each printer.}
\label{fig:gantt}
\end{figure}

Figure~\ref{fig:gantt} shows the schedules of the points with 3 and 2 jobs. With 3 jobs, printer 2 prints part 8, which has the latest due date, in a job of its own that starts after an idle time of 0.04 h and ends on its due date. With 2 jobs, part 8 joins the job of printer 1 and part 6 moves to the job of printer 2, so part 8 is finished 3.57 h early, part 6 is finished 1.55 h sooner, and the other parts finish at most 0.52 h later. All parts lie flat in both schedules, and the job heights are 44 to 52 mm. Every NSGA-II run found the points with 3 and 2 jobs, and the runs differed only in the points with one job. The average ratio of the hypervolume of an NSGA-II front to that of the exact front, defined in Section~\ref{sec:compare}, is 0.985, with a minimum of 0.970.

\subsection{Value of orientation}\label{sec:orientation}

We now compare the model, which chooses the orientations, with three fixed orientation policies, by minimizing $f_1$ alone, that is, by solving $P(\infty)$. The given policy keeps the orientation of the data ($o=1$). The laying policy puts every part in its flattest admissible orientation, and the standing policy puts it in its tallest one, as in \cite{oh2020impact}. Each policy is the model with every $\mathcal{O}_i$ reduced to one orientation, so its optimal $f_1$ is never below that of the model. We solve 10 instances for each combination of $n\in\{6,8,10\}$ parts and $|\mathcal{M}|\in\{1,2\}$ printers with the MILP, and 10 instances with 15, 20, and 30 parts and two printers with NSGA-II, taking the best of three runs for each policy. The solution of every policy is feasible for the model, so the column ``Optimized'' of Table~\ref{tab:policies} gives the best value found by the four runs of an instance.

\begin{table}[H]
\centering
\caption{Mean $f_1$ (h) over 10 instances under each orientation policy, with the increase over the optimized orientation (\%) in parentheses. The first six rows are solved by the MILP, and the last two columns give the number of instances in which the model with orientation choice was solved to optimality within 600 s and its mean solve time. The last three rows are solved by NSGA-II.}
\label{tab:policies}
\small
\begin{tabular}{rrrrrrrr}
\toprule
Printers & Parts & Optimized & Given & Laying & Standing & Optimal & Time (s) \\
\midrule
1 & 6 & 9.44 & 21.00 (122.4) & 9.44 (0.0) & 33.60 (255.8) & 10 & 4 \\
1 & 8 & 12.92 & 29.57 (128.9) & 13.18 (2.1) & 38.62 (199.0) & 10 & 46 \\
1 & 10 & 16.36 & 33.80 (106.6) & 16.36 (0.0) & 51.80 (216.6) & 7 & 363 \\
2 & 6 & 11.77 & 34.19 (190.5) & 11.77 (0.0) & 50.61 (330.0) & 10 & 10 \\
2 & 8 & 9.01 & 34.53 (283.2) & 9.01 (0.0) & 54.07 (500.0) & 10 & 57 \\
2 & 10 & 8.54 & 28.73 (236.5) & 8.54 (0.0) & 55.54 (550.5) & 2 & 538 \\
\midrule
2 & 15 & 11.74 & 34.55 (194.3) & 12.13 (3.3) & 63.69 (442.6) & & \\
2 & 20 & 18.96 & 37.73 (99.1) & 19.27 (1.7) & 51.83 (173.4) & & \\
2 & 30 & 31.45 & 52.06 (65.5) & 33.54 (6.7) & 63.88 (103.1) & & \\
\bottomrule
\end{tabular}
\end{table}

The orientation has a large effect. Over the 60 instances solved by the MILP, the given orientation raises the mean $f_1$ by 167\%, and standing the parts upright raises it by 318\%, because tall parts make long jobs. Put differently, choosing the orientation lowers the earliness and tardiness by 63\% compared with the given orientation and by 76\% compared with standing parts. Laying the parts flat is optimal in 58 of the 60 instances. In the other two, both with 8 parts and one printer, the model stands some parts upright so that more parts fit on the plate, and the laying policy is up to 8.5\% worse. The plate binds more often when there are more parts. With 15, 20, and 30 parts, the laying policy is worse than the best orientation found in 4, 3, and 7 of the 10 instances, by 3.3\%, 1.7\%, and 6.7\% on average and by up to 20.3\%. NSGA-II with orientation choice came within 0.4\% of these best values on average.

The model with orientation choice solves all instances with up to 8 parts, in at most 139 s, but only 9 of the 20 instances with 10 parts within 600 s. The fixed policies have a third of the binary variables and solve every instance. In the 11 unsolved instances, the best solution found by the solver is never better than the optimal solution of the laying policy and is worse in 4 of them, while the lower bound of the solver stays 35\% to 100\% below the best value. This weak bound comes from the big-M constraints \eqref{eqn:early} and \eqref{eqn:tardy}, and it is the reason why we use NSGA-II for larger instances.

\subsection{Sensitivity analysis}\label{sec:sensitivity}

\begin{figure}[H]
\centering
\includegraphics[width=\textwidth]{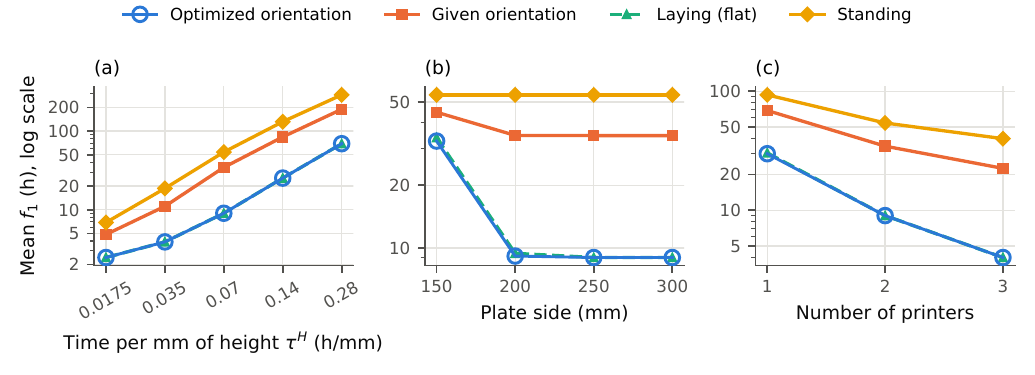}
\caption{Mean $f_1$ of the four orientation policies over 10 instances with 8 parts, when the layer time (a), the plate side (b), and the number of printers (c) change one at a time. The base values are $\tau^H=0.07$ h/mm, a 250 mm plate, and 2 printers. The dashed laying line lies under the optimized line where the two coincide.}
\label{fig:sens}
\end{figure}

We vary one parameter at a time for the 10 instances with 8 parts and 2 printers, keeping the parts and the due dates fixed: the layer time $\tau^H$ from a quarter to four times its base value, the side of the plate from 150 mm to 300 mm (with 80\% usable area), and the number of printers from 1 to 3. All 400 MILP solves of this section, including the 100 with orientation choice, were optimal, in at most 350 s. Figure~\ref{fig:sens} shows the mean $f_1$ of each policy.

The layer time scales the cost of tall jobs. When $\tau^H$ grows from 0.0175 to 0.28 h/mm, the mean $f_1$ with optimized orientation grows from 2.46 h to 69.13 h, and the model uses fewer and fuller jobs, 2.4 on average instead of 4.5, because every job adds the layer time of its tallest part. The given orientation is 96\% worse at the smallest layer time and 283\% worse at the base value, so the orientation matters more when the layers are slow. In absolute terms, the loss of the given orientation grows from 2.4 h to 120 h. The laying policy stays optimal at every layer time, since the plate is large enough.

The plate size decides whether laying parts flat is enough. With a 250 mm or 300 mm plate, laying is optimal in every instance, and the larger plate does not change $f_1$. With a 200 mm plate, laying is worse in 2 of the 10 instances, by up to 21.5\%. With a 150 mm plate, many flat orientations no longer fit, and $f_1$ rises to 32.51 h. Laying is then worse in 8 of the 10 instances, by 3.8\% on average and by up to 12.6\%, because the model stands some parts upright and needs 3.3 jobs on average instead of 3.7. The benefit of a larger printer therefore depends on the size of the parts. Beyond the size at which the flat parts fit, a larger plate gives no gain.

More printers reduce the earliness and tardiness sharply when the due dates stay the same. The mean $f_1$ falls from 29.83 h with one printer to 9.01 h with two and 4.00 h with three, while the gain of the second printer is much larger than that of the third. With one printer, laying is worse in 2 of the 10 instances, by up to 11.2\%. Across all settings, the given and standing orientations are 37\% to 461\% and 66\% to 897\% worse than the optimized orientation.

\subsection{NSGA-II and the epsilon-constraint method}\label{sec:compare}

We compare NSGA-II with the epsilon-constraint method on 5 instances for each number of parts $n\in\{8,10,15,20,30,40\}$, all with 2 printers. The grid version of the epsilon-constraint method has a time limit of 60 s for each of its 22 stage solves. For the instances with 8 parts, the complete version, with 600 s per solve, gives the whole Pareto front. NSGA-II is run 10 times on each instance with different random seeds.

Since the Pareto front of a larger instance is unknown, every front is compared with the reference set $\mathcal{R}$ of its instance, which consists of the nondominated points of the union of all fronts found for the instance. Both objectives are scaled to $[0,1]$ by their smallest and largest values in $\mathcal{R}$. The hypervolume of a front is the area of the scaled objective space that is dominated by the front and bounded by the point $(1.1,1.1)$ \cite{zitzler2003performance}. We report it as a ratio to the hypervolume of $\mathcal{R}$, so a front with ratio 1 dominates the same area as $\mathcal{R}$. The inverted generational distance (IGD) is the mean Euclidean distance from a point of $\mathcal{R}$ to the nearest point of the front in the scaled space \cite{coello2005solving}. A small IGD means that every point of $\mathcal{R}$ has a point of the front nearby. Table~\ref{tab:compare} gives the averages over the instances of each size.

\begin{table}[H]
\centering
\caption{Epsilon-constraint grid and NSGA-II on 5 instances of each size with 2 printers. The hypervolume is the ratio to that of the reference set, given for NSGA-II as the mean over the 50 runs with the standard deviation in parentheses. The column ``Optimal'' gives the share of the stage solves of the grid that were optimal within 60 s. The times are per instance for the grid and per run for NSGA-II, and ``Points'' is the number of nondominated points returned.}
\label{tab:compare}
\small
\setlength{\tabcolsep}{3.5pt}
\begin{tabular}{rrrrrrrrrr}
\toprule
 & \multicolumn{5}{c}{Epsilon-constraint grid} & \multicolumn{4}{c}{NSGA-II} \\
\cmidrule(lr){2-6}\cmidrule(lr){7-10}
Parts & Hypervolume & IGD & Optimal (\%) & Time (s) & Points & Hypervolume & IGD & Time (s) & Points \\
\midrule
8 & 1.000 & 0.030 & 98 & 93 & 3.8 & 0.992 (0.011) & 0.051 & 38 & 3.8 \\
10 & 0.995 & 0.083 & 90 & 248 & 4.0 & 0.955 (0.051) & 0.073 & 40 & 5.9 \\
15 & 0.965 & 0.185 & 35 & 984 & 4.6 & 0.992 (0.010) & 0.099 & 43 & 8.8 \\
20 & 0.884 & 0.156 & 4 & 1323 & 5.2 & 0.986 (0.020) & 0.074 & 45 & 12.8 \\
30 & 0.490 & 0.574 & 0 & 1381 & 4.0 & 0.989 (0.010) & 0.068 & 53 & 19.2 \\
40 & 0.523 & 0.503 & 0 & 1381 & 4.0 & 0.976 (0.016) & 0.073 & 53 & 20.8 \\
\bottomrule
\end{tabular}
\end{table}

With 8 parts, the complete method found the whole front, 5.0 points on average, in 72 s on average, and every solve was optimal. The grid missed 6 of the 25 Pareto points, all with one job, and kept at least 0.998 of the hypervolume. An NSGA-II run found 3.1 of the 5.0 Pareto points on average, and the 10 runs together found 22 of the 25 points. The 3 points that no run found also have one job. The hypervolume ratio of NSGA-II is 0.992 on average and at least 0.967, and in every run the smallest $f_1$ of the front is the optimal value of $P(\infty)$ from Section~\ref{sec:orientation}.

With 10 parts, 90\% of the stage solves of the grid were optimal, and the grid has a higher hypervolume ratio than NSGA-II, 0.995 against 0.955. The difference lies in the points with one job. When all parts are in one job, minimizing $f_2$ is a multiple-choice subset sum problem over the orientations, which the MILP solves exactly for $\varepsilon$ close to 0, while NSGA-II finds the best of these points only in some runs. Since the points with one job cover the upper part of the scaled $f_1$ range, a run that misses them loses much hypervolume, and the worst run keeps 0.810. NSGA-II, in turn, has the lower IGD, because the first stage of the grid, which minimizes $f_1$ alone, stopped at 60 s in 4 of the 5 instances with a value 7\% to 75\% above the smallest $f_1$ of NSGA-II. All 10 NSGA-II runs reach this smallest value on each instance. It is the optimal value of $P(\infty)$ for one instance, and for the other four it equals the best value found in Section~\ref{sec:orientation}, which the MILP with orientation choice did not reach within 600 s in two of them.

\begin{figure}[H]
\centering
\includegraphics[width=\textwidth]{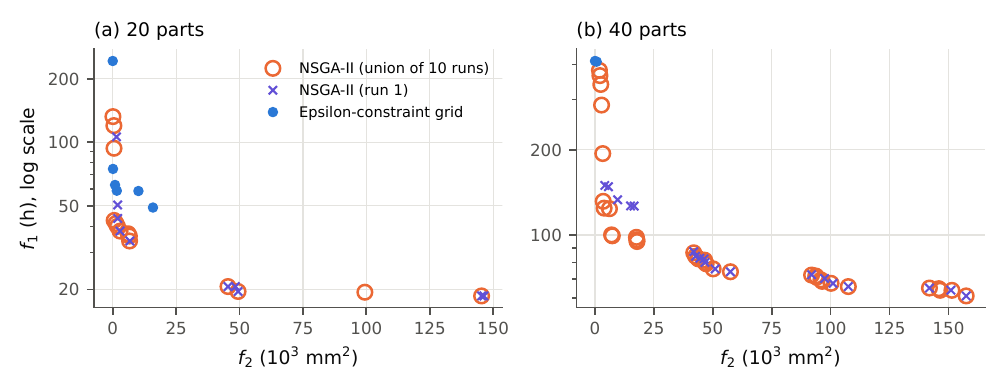}
\caption{Fronts of instance 0 with 20 parts (a) and 40 parts (b). The circles are the nondominated points of the union of the 10 NSGA-II runs, the crosses are the front of a single run, and the dots are the front of the epsilon-constraint grid.}
\label{fig:fronts}
\end{figure}

From 15 parts on, the grid falls behind. The share of optimal stage solves drops to 35\% with 15 parts, 4\% with 20 parts, and none with 30 and 40 parts, and the hypervolume ratio of the grid falls to about one half with 30 and 40 parts, after 1381 s per instance. The smallest $f_1$ of the grid is 1.4 to 3.1 times that of NSGA-II with 15 and 20 parts and 3.4 to 9.8 times with 30 and 40 parts. NSGA-II keeps a mean hypervolume ratio of 0.976 to 0.992 and an IGD below 0.1, with 9 to 21 points per front, in 43 s to 53 s per run on average, and the union of its runs supplies 89\% to 92\% of the points of the reference sets with 20 or more parts. Figure~\ref{fig:fronts} shows the fronts of instance 0 with 20 and 40 parts. The grid finds points only at small $f_2$, where the jobs are nearly full, and its smallest $f_1$ is 2.6 and 6.7 times that of NSGA-II on these two instances. The fronts of NSGA-II also cover the plans with spare plate area, which have the smallest $f_1$.

\section{Conclusions}\label{sec:conclusion}

This paper studied the allocation of parts to jobs on identical powder bed 3D printers, together with the orientation of each part and the sequence of the jobs, with two objectives, the weighted earliness and tardiness of the parts and the unused plate area. We formulated a bi-objective MILP model in which the orientation enters through the assignment variables, so that no product of variables has to be linearized. We showed that the unused area equals the plate area times the number of jobs minus the total base area, gave a bound on completion times that keeps an optimal solution and tightens the big-M constraints, and showed that the optimal timing of a fixed job sequence is an isotonic regression solved exactly by PAVA. The epsilon-constraint method gives the whole Pareto front of small instances, and NSGA-II, whose decoder uses the timing result, handles larger ones.

The experiments show that the orientation is a first-order decision. On average, keeping the orientation of the data raises the earliness and tardiness by 167\%, and standing the parts upright raises it by 318\%. Laying the parts flat is optimal when the plate is large compared with the parts. When the plate binds, because it is small, because there is one printer, or because there are many parts, the optimal plan stands some parts upright, and laying all parts flat is up to 21.5\% worse. The Pareto front is organized by the number of jobs. For a given number of jobs, filling the plate by reorienting parts costs much time for little area, so the useful plans are the best ones for each number of jobs. A slower layer time raises the value of the orientation, a larger plate helps only up to the size at which the flat parts fit, and a second printer reduces the earliness and tardiness much more than a third one.

NSGA-II with the timing decoder reaches the smallest known earliness and tardiness in every run on the instances with 8 and 10 parts. From 15 to 40 parts, its fronts keep 98\% to 99\% of the hypervolume of the best known fronts in less than a minute per run on average, while the epsilon-constraint method with a time limit keeps about half of it with 30 and 40 parts after 23 minutes. The exact method remains the better choice for small instances and for the plans with one job.

\section*{Declaration of competing interest}

The author declares no known competing financial interests or personal relationships that could have appeared to influence the work reported in this paper.

\section*{Data and code availability}

The code and data that produce the results of Section~\ref{sec:numexp} are provided as supplementary material.

\bibliography{references}

\end{document}